\documentclass[12pt]{elsarticle}

\usepackage[english]{babel}
\usepackage{csquotes}
\usepackage{hyperref}

\usepackage[T1]{fontenc}

\usepackage{amssymb}
\usepackage{dcolumn}
\usepackage{enumerate}
\usepackage{mathdots}
\usepackage{tabularx}
\usepackage{subfigure}
\usepackage{amsmath}
\usepackage{amsthm}
\usepackage{tikz} 
\usepackage{float}
\usepackage{mathtools}
\usepackage{graphicx}
\usepackage{pgfplots}
\usepackage{natbib}

\usepackage[margin=3.5cm]{geometry}

\usepackage{pgf,tikz}
\usepackage{mathrsfs}

\usepackage{algorithm,algpseudocode,float}
\allowdisplaybreaks

\newcommand{\bn}{B_n}

\newcommand{\pn}{\mathbb{P}_n}
\newcommand{\bpn}{\mathbb{BP}_n}

\newcommand{\Sp}{\text{sp}}

\newcommand{\bp}[1]{\mathbb{BP}_{#1}}
\newcommand{\rvline}{\hspace*{-\arraycolsep}\vline\hspace*{-\arraycolsep}}
\newcommand{\ol}[1]{\overline{#1}}

\newcommand{\indsize}{\scriptsize}
\newcommand{\colind}[2]{\displaystyle\smash{\mathop{#1}^{\raisebox{.5\normalbaselineskip}{\indsize #2}}}}
\newcommand{\rowind}[1]{\mbox{\indsize #1}}

\usepackage{pgf,tikz}\usepackage{mathrsfs}\usetikzlibrary{arrows}

\usetikzlibrary{matrix,fit, arrows, shapes}

\numberwithin{equation}{section}

\newtheorem{thm}{Theorem}
\newtheorem{prop}[thm]{Proposition}
\newtheorem{cor}[thm]{Corollary}
\newtheorem{lem}[thm]{Lemma}

\newtheorem{conj}[thm]{Conjecture}

\theoremstyle{definition}
\newtheorem{defn}{Definition}[section]

\pgfplotsset{compat=1.18}

\begin{document}
 
\title{Some integer values in the spectra of burnt pancake graphs}

\author[1]{Sa\'ul A. Blanco}
\affiliation[1]{organization={Department of Computer Science, Indiana University}, city={Bloomington, IN},
postcode={47408}, country={U.S.A.}}
\ead{sblancor@iu.edu}

\author[2]{Charles Buehrle}
\affiliation[2]{organization={Department of Mathematics, Physics, and Computer Studies, Notre Dame of Maryland University}, city={Baltimore, MD}, postcode={21210}, county={U.S.A.}}
\ead{cbuehrle@ndm.edu}

\date{September 3, 2024}

\begin{abstract} 
The burnt pancake graph, denoted by $\bpn$, is formed by connecting signed permutations via prefix reversals. Here, we discuss some spectral properties of $\bpn$. More precisely, we prove that the adjacency spectrum of $\bpn$ contains all integer values in the set $\{0, 1, \ldots, n\}\setminus\{\left\lfloor n/2 \right\rfloor\}$.

\end{abstract}

\maketitle

\section{Introduction}

Let $n$ be a positive integer and $[n]$ denote the set $\{1,2,\ldots,n\}$. Similarly, for integers $a$ and $b$, with $a<b$, $[a,b]$ denotes the set $\{a,a+1,\ldots, b\}$. In addition, let $S_n$ denote the symmetric group of order $n!$ and $B_n$ denote the \emph{hyperoctahedral group} of order $2^nn!$. In other words, if $\pm[n]$ denotes the set $\{-n,-(n-1),\ldots,-1,1,2,\ldots,n\}$, then $B_n$ is the group of all bijections $\sigma$ on $\pm[n]$ satisfying $\sigma(-i)=-\sigma(i)$ for $1\leq i\leq n$, with the group operation being composition. We shall utilize subscripts for the image of an element in a permutation, that is $\sigma_i = \sigma(i)$ for $i\in[n]$, and denote negative signs with a bar on top of the character, for example, $\overline{3}=-3$. If $\pi=\pi_1\pi_2\cdots\pi_n$ denotes a permutation of length $n$ in one-line notation, then the $i$-th prefix reversal $r_i:S_n\to S_n$ is defined by $r_i(\pi)=\pi_i\pi_{i-1}\cdots\pi_1\pi_{i+1}\cdots\pi_n$, for $2\leq i\leq n$. Similarly, we also write $\sigma\in B_n$ in one-line notation as $\sigma=\sigma_1\sigma_2\cdots\sigma_n$,  for example, a particular element of $B_5$ is $3\overline{2}4\overline{5}1$. The $i$-th \emph{signed prefix-reversal} $r_i$, with $1\leq i\leq n$, is defined as follows: $r_i:B_n\to B_n$ is given by $r_i(\sigma_1\cdots\sigma_n)=\ol{\sigma_i}\, \ol{\sigma_{i-1}}\cdots\ol{\sigma_1}\sigma_{i+1}\cdots\sigma_n$. The prefix reversals $r_i$ defined above are themselves elements of $S_n$ (or in the case of signed reversals, $B_n$). Thus, they may be expressed in cyclic notation with $r_i = (1,i)(2,i-1)\cdots(\left\lfloor n/2 \right\rfloor, \left\lceil n/2 \right\rceil) \in S_n$ for $2 \leq i \leq n$ and $r_i = (1,\ol{i})(2, \ol{i-1}) \cdots (i,\ol{1})\in\bn$ for $1 \leq i \leq n$. 

The \emph{pancake graph}, denoted by $\pn$, has the elements of $S_n$ as its set of vertices and its edge set is $\{(\pi,r_i(\pi)):\pi\in S_n, 2\leq i\leq n\}$. In addition, the \emph{burnt pancake graph}, denoted by $\bpn$, has the elements of $B_n$ as its vertex set and its edge set is $\{(\sigma,r_i(\sigma):\sigma\in B_n, 1\leq i\leq n\}$.

There is a great deal that is known about the pancake graph $\pn$, including its cyclic structure, see Kanevsky and Feng as well as the work of Konstantinova and Medvedev~\cite{KF95,KM10,KM11, KM16}. Moreover, the existence of some integer eigenvalues in $\pn$ was established in Dalf\'{o} and Fiol~\cite{Dalfo}. More precisely, they establish the following proposition.

\begin{prop}(~\cite[Proposition 2.2]{Dalfo})
    The spectrum of $\pn$ with $n\geq3$ contains every element in the set $[-1,n-1]\setminus\{\lfloor (n-2)/2 \rfloor\}$.
\end{prop} If the case is $n=1$ and $n=2$, then $\pn$ is isomorphic to a vertex and an edge, respectively. 

We recall that if $G$ is a group with generating set $S$, then the \emph{Cayley graph} of $G$ with respect to $S$, denoted by $Cay(G,S)$ is the graph with vertex set $G$ and edge set $\{(g,sg):g\in G, s\in S\}$. Notice that $\pn$ is the Cayley graph of $S_n$ with respect to the generators $\{r_i\}_{i=2}^n$ and $\bpn$ is the Cayley graph of $B_n$ with respect to $\{r_i\}_{i=1}^n$. 

Both $\pn$ and $\bpn$ have been widely studied objects due to their connections to parallel computing, see Kanevsky and Feng or Lakshmivarahan, Jwo, and Dhall~\cite{KF95, LJD93}, and bioinformatics, see Fertin, Labarre, Rusu, Tannier, and Vialette or Hannenhalli and Pevzner~\cite{FLRTV, HannenPev}. In particular, the cyclic structure of $\bpn$, and in general of \emph{prefix reversal graphs}, is studied in the authors' work, also with Patidar,~\cite{BB23, BBP19, BBP19Perm}.

\subsection*{Our contribution.} In this paper, we show that the burnt pancake graph $\bpn$ has all integers in $[0,n]\setminus\{\lfloor n/2\rfloor\}$ as eigenvalues of its adjacency matrix.

\subsection{Regular partitions} Given a graph $G=(V,E)$, we define the \emph{neighborhood} of $u\in V$ as the set $N(u)=\{v\in V: (u,v)\in E\}$. In other words, the set of all vertices adjacent to vertex $u$. We follow the notation from Dalf\'{o} and Fiol~\cite{Dalfo}. For a graph $G$, let $A_G$ denote the adjacency matrix of $G$. The set 
\[
\Sp(G)=\Sp(A_G)=\{[\lambda_0]^{m_0},[\lambda_1]^{m_1},\ldots,[\lambda_d]^{m_d}\},
\] where the set $\{\lambda_i\}_{i=1}^d$ is the set of distinct eigenvalues of $A_G$ and $m_i$ is the multiplicity of $\lambda_i$ for $1\leq i\leq d$, is called the \emph{(adjacency) spectrum} of $G$, or of $A_G$.

Given a graph $G=(V,E)$, by a \emph{partition} $P$ of $V$, we mean the subsets $P=\{V_1,\ldots, V_k\}$ of $V$ satisfying $V_1\cup\cdots\cup V_k=V$ and $V_i\cap V_j=\emptyset$ if $i\neq j$. Furthermore, we say that $P$ is \emph{equitable} (also referred to as \emph{regular}) if, for every $w_{i,1}, w_{i,2}\in V_i$, $|N(w_{i,1})\cap V_j|=|N(w_{i,2})\cap V_j|$ for $1\leq i,j\leq k$. Moreover, let $\mathbf{s}_{u}$ be the $k\times 1$ characteristic column vector corresponding to $u\in V$; namely for $1\leq i\leq k$, $s_{i,u}=1$ if $u\in V_i$ and $s_{i,u}=0$ if $u\not\in V_i$. The matrix $S_{P}=(\mathbf{s}_{i,u})_{k\times |V|}$ is called the \emph{characteristic matrix} of $P$.

Let us assume the partition $P=\{V_1, \ldots, V_k\}$ is an equitable partition of $V$. We let $b_{i,j}=|N(u)\cap V_j|$, where $u\in V_i$ and $B=(b_{i,j})_{k\times k}$ is referred to as the \emph{quotient matrix} of $A_G$ with respect to $P$.

The quotient matrix $B$ can be used to derive information about the spectrum of the graph $G$. Indeed, we summarize the necessary results that we utilize here.

\begin{lem}[Lemma 1.1 in~\cite{Dalfo}]\label{lem:quotient}
    Let $G=(V,E)$ be a graph with adjacency matrix $A_G$ and $P$ be a partition of $V$ with characteristic matrix $S_{P}$. Then the following statements are true.
    
    (i) The partition $P$ is regular if and only if there exists a matrix $C$ such that $S_{P}C=A_GS_{P}$. Moreover, $C=B$, the quotient matrix of $A_G$ with respect to $P$.

    (ii) If $P$ is a regular partition and $\mathbf{v}$ is an eigenvector of $B$, then $S_{P}\mathbf{v}$ is an eigenvector of $A_G$. In particular, $\Sp(B)\subseteq\Sp(A_G)$.
\end{lem}

Furthermore, the following proposition is used.

\begin{prop}[Proposition 2.1 in~\cite{Dalfo}]\label{prop:subgroupofsym}
    Let $G=\text{Cay}(Gr,S)$ with $Gr$ being a subset of the symmetric group and the generating set $S$ being a set of permutations $\{\sigma_1,\sigma_2,\ldots,\sigma_\ell\}$. Then, $G$ has a regular partition with quotient matrix $B=\sum_{i=1}^{\ell}P(\sigma_i)$, where $P(\sigma)$ denotes the permutation matrix of $\sigma$.
\end{prop}

\section{Some integer eigenvalues}

With the previous results and notation established, in this section, we show that $\bpn$ has all integer eigenvalues in the set $[0,n]\setminus\{\lfloor{n}/{2}\rfloor\}$.

Notice that $\bpn$ can be naturally thought of as a subgroup of $S_{2n}$ with generating set $\{r_i\}_{i=1}^n$. Thus, we can apply Proposition~\ref{prop:subgroupofsym} to obtain a matrix corresponding to an equitable partition of $\bpn$. In cyclic notation, \[
r_i=(1,\overline{i})(2,\overline{i-1})\cdots (i,\overline{1}).
\]

For an integer $k$, the \emph{$k$-th diagonal} of a matrix refers to the set of entries that lie on a diagonal that is $k$ positions off from the main diagonal. So, the main diagonal corresponds to $k=0$, $k>0$ corresponds to a diagonal $k$ positions above the main diagonal, and $k<0$ corresponds to a diagonal that is $|k|$ positions below the main diagonal.

Ordering the elements of $\pm[n]=\{\overline{n},\overline{n-1},\ldots,\overline{1},1,2,\ldots,n\}$  increasingly and indexing the permutation matrix $P(r_i)$ of $r_i$ accordingly, $P(r_i)$ has the following description. Since $r_i$ leaves fixed any element $j$ such that $|j|>i$, $P(r_i)$ has a 1 in the diagonal entry for the first and last $n-i$ rows. Then, there is a 1 in the $\overline{i}$-th diagonal entries indexed from $(\overline{i},1)$ to $(\overline{1},i)$ and, symmetrically, there is a 1 in the $i$-th diagonal entries indexed from $(1,\overline{i})$ to $(i,\overline{1})$. Any other entry in $P(r_i)$ is 0. We give an example of these permutation matrices in the case $\bp{3}$.

\[
  P(r_1) =
  \begin{array}{@{}c@{}}
    \rowind{$\overline{3}$} \\ \rowind{$\overline{2}$} \\ \rowind{$\overline{1}$} \\ \rowind{1} \\ \rowind{2} \\\rowind{3}
  \end{array}
  \left[
  \begin{array}{ *{6}{c} }
     \colind{1}{$\overline{3}$}  &  \colind{0}{$\overline{2}$}  &  \colind{0}{$\overline{1}$}  & \colind{0}{1} & \colind{0}{2} & \colind{0}{3} \\
0 & 1 & 0  & 0 & 0 & 0 \\
 0 & 0 & 0  & 1 & 0 & 0 \\
 0 & 0 & 1  & 0 & 0 & 0 \\
 0 & 0 & 0  & 0 & 1 & 0 \\
 0 & 0 & 0  & 0 & 0 & 1
  \end{array}
\right],\quad P(r_2)=
\begin{array}{@{}c@{}}
    \rowind{$\overline{3}$} \\ \rowind{$\overline{2}$} \\ \rowind{$\overline{1}$} \\ \rowind{1} \\ \rowind{2} \\\rowind{3}
  \end{array}
  \left[
  \begin{array}{ *{6}{c} }
     \colind{1}{$\overline{3}$}  &  \colind{0}{$\overline{2}$}  &  \colind{0}{$\overline{1}$}  & \colind{0}{1} & \colind{0}{2} & \colind{0}{3} \\
 0 & 0 & 0  & 1 & 0 & 0 \\
 0 & 0 & 0  & 0 & 1 & 0 \\
 0 & 1 & 0  & 0 & 0 & 0 \\
 0 & 0 & 1  & 0 & 0 & 0 \\
 0 & 0 & 0  & 0 & 0 & 1 \\
\end{array}
\right], \text{ and }
\] \vspace{0.2cm}
\[
P(r_3)=\begin{array}{@{}c@{}}
    \rowind{$\overline{3}$} \\ \rowind{$\overline{2}$} \\ \rowind{$\overline{1}$} \\ \rowind{1} \\ \rowind{2} \\\rowind{3}
  \end{array}
  \left[
  \begin{array}{ *{6}{c} }
     \colind{0}{$\overline{3}$}  &  \colind{0}{$\overline{2}$}  &  \colind{0}{$\overline{1}$}  & \colind{1}{1} & \colind{0}{2} & \colind{0}{3} \\
 0 & 0 & 0  & 0 & 1 & 0 \\
 0 & 0 & 0  & 0 & 0 & 1 \\
 1 & 0 & 0  & 0 & 0 & 0 \\
 0 & 1 & 0  & 0 & 0 & 0 \\
 0 & 0 & 1  & 0 & 0 & 0 \\
\end{array}
\right].
\]

With this clear pattern of the permutation matrices, we can describe the matrix $B=\sum_{i=1}^nP(r_i)$ from Proposition~\ref{prop:subgroupofsym}.

Let $A_n$ and $C_n$ be the $n\times n$ diagonal matrices with diagonal entries $a_{i,i}=n-i$ and $c_{i,i}=i-1$ for $1\leq i\leq n$, respectively. Furthermore, let $D_n$ be the upper triangular matrix where every entry in the diagonal and above the diagonal is 1. We may then define the $2n\times 2n$ block matrix $M(\bpn)$ as follows:

\[
M(\bpn)=\begin{bmatrix}A_n& D^T_n\\ D_n& C_n\end{bmatrix},
\]
where, of course, $D^T_n$ is the transpose of $D_n$.

For example, the special case where $n=3$ is shown below. The lines are included to better identify the blocks $A_n$, $C_n$, and $D_n$, and we have

\[
    M(\bp{3})= \left[ \begin{matrix}
        \begin{matrix}
            2 & 0 & 0\\
            0 & 1 & 0\\
            0 & 0 & 0
        \end{matrix} & \rvline &
        \begin{matrix}
            1 & 0 & 0\\
            1 & 1 & 0\\
            1 & 1 & 1 
        \end{matrix} \\ \hline
        \begin{matrix}
            1 & 1 & 1\\
            0 & 1 & 1\\
            0 & 0 & 1
        \end{matrix} & \rvline &
        \begin{matrix}
            0 & 0 & 0\\
            0 & 1 & 0\\
            0 & 0 & 2
        \end{matrix}
    \end{matrix} \right]. 
\]
One readily verifies that $M(\bp{3})=P(r_1)+P(r_2)+P(r_3)$. In general, we have the following lemma, which shows that the quotient matrix of Proposition~\ref{prop:subgroupofsym}, for the generating set of prefix reversals, is precisely this block matrix, $M(\bpn)$.

\begin{lem} For $n\geq1$,
    $M(\bpn)=\sum_{i=1}^nP(r_i)$. 
\end{lem} 

\begin{proof} Let $M_n=\sum_{i=1}^nP(r_i)$.  When indexing the entries of $M_n$, we take the standard ordering on the elements of $\pm[n]$, namely $\ol{n}<\ol{n-1}<\cdots<\ol{1}<1<\cdots <n$ (see, for example, $P(r_1)$ earlier in this section). Let us recall that $r_i=(1,\overline{i})(2,\overline{i-1})\cdots (i,\overline{1})$, for $1\leq i\leq n$, and let $R_n=\{r_i:1\leq i\leq n\}$. We will argue by cases depending on the indices of the entries in $M_n$.

\textbf{Diagonal entries.} From the cyclic-notation representation of $r_i$, we see that $r_i$ leaves every $j\in\pm[n]$ with $|j|>i$ fixed. Therefore, it follows that there are $|j|-1$ elements from $R_n$ that leave $|j|$ fixed, corresponding to the diagonal entries of $M_n$. Thus the diagonal elements of $M_n$ are given by the diagonal elements of $A_n$ and $C_n$ seen as block matrices in $M(\bpn)$. So for diagonal entries, $M_n$ and $M(\bpn)$ coincide. 

\textbf{Non-diagonal entries indexed by two positive or two negative elements in $\pm[n]$.} Notice that $r_i$ swaps $j$ and $\overline{i-j+1}$, with $1\leq j\leq i$. Thus the non-diagonal entries of $P(r_i)$ are indexed by pairs $(a,b)$ where exactly one of $a$ or $b$ is positive and the other one is negative. Hence, every non-diagonal entry of $M_n$ indexed by two positive or two negative elements must necessarily be 0, thus coinciding with $A_n$ and $C_n$ seen as block matrices in $M(\bpn)$. 

\textbf{Non-diagonal entries indexed by one positive and one negative element in $\pm[n]$.} Notice that the non-zero, non-diagonal entries of $P(r_i)$ lie on the $i$-th and $\overline{i}$-th diagonal. In fact, there are $i$ entries of 1 in the $\overline{i}$-th diagonal from indices $(1,\overline{i})$ to $(i,\overline{1})$ and there are $i$ entries of 1 in the $i$-th diagonal from indices $(\overline{i},1)$ to $(\overline{1},i)$, respectively. Moreover, notice that for any prefix reversals, $r_{i_1}$ and $r_{i_2}$, each swap different elements in $\pm[n]$ if $i_1\neq i_2$. So, the entries in $M_n$ that are indexed by one positive and one negative element are given by $D_n$ and $D^T_n$ seen as block matrices in $M(\bpn)$. 

Therefore, $M_n=M(\bpn)$, as claimed.    \end{proof}

In particular, from Lemma~\ref{lem:quotient}(ii), and Proposition~\ref{prop:subgroupofsym}, $\Sp(M(\bpn))\subseteq \Sp(\bpn)$. We are now ready to state and prove our main result. 

\begin{thm}
    The spectrum of $\bpn$ includes $[0,n]\setminus\{\lfloor {n}/{2}\rfloor\}$.
\end{thm}

\begin{proof}
    One can directly verify that the following are eigenvalues and their corresponding eigenvectors of $M(\bpn)$. 
    For even $n$,
\begin{align*}
    n&\text{ and }(\underbrace{1,1,\ldots,1}_{2n}),\\
    n-1&\text{ and }(n-2,\underbrace{-1,\ldots,-1}_{n-2},0,0,\underbrace{-1,\ldots,-1}_{n-2},n-2),\\
    \vdots &\hspace{2cm}\vdots \\
     n-i&\text{ and }(\underbrace{0,\ldots,0}_{i-1}, n-2i,\underbrace{-1,\ldots,-1}_{n-2i},\underbrace{0,\ldots,0}_{2i},\underbrace{-1,\ldots,-1}_{n-2i},n-i,\underbrace{0,\ldots,0}_{i-1})\\
     &\hspace{-60pt}\text{ for $i$ with }1\leq i\leq  {n}/{2} -1,\\ 
     \frac{n}{2}-1&\text{ and }(\underbrace{0,\ldots,0}_{\frac{n}{2}-1},1,-1\underbrace{0,\ldots,0}_{2(\frac{n}{2}-1)},-1,1,\underbrace{0,\ldots,0}_{\frac{n}{2}-1}),\\
          \frac{n}{2}-1-j&\text{ and }(\underbrace{0,\ldots,0}_{\frac{n}{2}-1-j},\underbrace{1,\ldots,1}_{2j+1},-2j-1,\underbrace{0,\ldots,0}_{2(\frac{n}{2}-1-j)},-2j-1,\underbrace{1,\ldots,1}_{2j+1},\underbrace{0,\ldots,0}_{\frac{n}{2}-1-j})
\end{align*} for $j$ with $1\leq j\leq {n}/{2}-1$.

Furthermore, for $n$ odd, the eigenvalues and corresponding eigenvectors are as follows.
\begin{align*}
    n&\text{ and }(\underbrace{1,1,\ldots,1}_{2n}),\\
    n-1&\text{ and }(n-2,\underbrace{-1,\ldots,-1}_{n-2},0,0,\underbrace{-1,\ldots,-1}_{n-2},n-2),\\
    \vdots &\hspace{2cm}\vdots \\
     n-i&\text{ and }(\underbrace{0,\ldots,0}_{i-1}, n-2i,\underbrace{-1,\ldots,-1}_{n-2i},\underbrace{0,\ldots,0}_{2i},\underbrace{-1,\ldots,-1}_{n-2i},n-i,\underbrace{0,\ldots,0}_{i-1})\\
     &\hspace{-64pt}\text{ for $i$ with }1\leq i\leq \left\lfloor {n}/{2}\right\rfloor,\\ 
     \left\lfloor\frac{n}{2}\right\rfloor-1&\text{ and }(\underbrace{0,\ldots,0}_{\left\lfloor\frac{n}{2}\right\rfloor-1},1,-1\underbrace{0,\ldots,0}_{2(\left\lfloor\frac{n}{2}\right\rfloor-1)},-1,1,\underbrace{0,\ldots,0}_{\left\lfloor\frac{n}{2}\right\rfloor-1}),\\
          \left\lfloor\frac{n}{2}\right\rfloor-1-j&\text{ and }(\underbrace{0,\ldots,0}_{\left\lfloor\frac{n}{2}\right\rfloor-1-j},\underbrace{1,\ldots,1}_{2j+1},-2j-1,\underbrace{0,\ldots,0}_{2(\left\lfloor\frac{n}{2}\right\rfloor-1-j)},-2j-1,\underbrace{1,\ldots,1}_{2j+1},\underbrace{0,\ldots,0}_{\left\lfloor\frac{n}{2}\right\rfloor-1-j})
     \end{align*}  for $j$ with $1\leq j\leq \left\lfloor{n}/{2}\right\rfloor-1$.
\end{proof}

\section{Final remarks}

An intriguing open question regarding the burnt pancake graph is to determine its \emph{spectral gap}, defined as the difference between the largest two eigenvalues of the adjacency matrix of $\bpn$. This question, with some interesting numerical results, was also mentioned in Chung and Tobin~\cite{CT17}, where those authors utilize graph covering methods to establish that the spectral gap of a family of graphs that contains $\pn$ is 1. We investigated several graph coverings/projections, but $\bpn$ proved to be more idiosyncratic. Each potential projection (weighted graphs covered by $\bpn$ with particular conditions on weights) we attempted yielded integer eigenvalues, already established in our main result above. However, through a particular projection of $\bpn$, we are able to establish that the multiplicity of the eigenvalue $n-1$ is at least 2, which is to be expected. Here, we provide a brief overview of the particular graph covering/projection we investigated.

\begin{defn}\label{d:cover}
    We follow the notation of Chung and Tobin as well as Chung and Yau~\cite{CT17, CY99}. Let $G$ be a weighted graph. Then, we say that the weighted graph $\tilde{G}$ is a \emph{covering} of $G$, or that $G$ is a \emph{projection} of $\tilde{G}$, if there exists a surjection $p:V(\tilde{G})\to V(G)$ that satisfies 
\begin{enumerate}
    \item For all $\tilde{v}_1,\tilde{v}_2\in V(\tilde{G})$ with $p(\tilde{v}_1)=p(\tilde{v}_2)$ and for all $v\in V(G)$,
    \[
    \sum_{\tilde{v}_3\in p^{-1}(v)}w(\tilde{v}_1,\tilde{v}_3)=\sum_{\tilde{v}_3\in p^{-1}(v)}w(\tilde{v}_2,\tilde{v}_3). 
    \]
    \item There exists $m\in\mathbb{R}^+\cup\{\infty\}$ such that for all $v_1,v_2\in V(G)$,
    \[
    \sum_{\substack{u\in p^{-1}(v_1)\\ v\in p^{-1}(v_2)}}w(u,v)= mw(v_1,v_2).
    \] Here, $m$ is referred to as the \emph{index} of $p$.
\end{enumerate}
\end{defn}

The \emph{Laplacian} of a graph is the matrix $L = D - A_G$, where $\nu = |V(G)|$, $D$ is a $(\nu \times \nu)$-diagonal matrix whose entries are the degrees of the corresponding vertices and $A_G$ is the graph's adjacency matrix. Furthermore, the \emph{normalized Laplacian} is defined as $\mathcal{L} = D^{-1/2} L D^{-1/2}$. When dealing with a regular graph of degree $d$ the eigenvalues of $\mathcal{L}$, $0=\lambda_0 \leq \lambda_1 \leq \ldots \leq \lambda_{\nu-1}$, are related to the eigenvalues of the adjacency matrix $A_G$, $\mu_1 \geq \mu_2 \geq \ldots \geq \mu_{\nu}$, by $\mu_i = d(1-\lambda_i)$ for all $i$. A consequence of the Covering-Correspondence Theorem from Chung and Tobin~\cite[Theorem 3]{CT17} is the following result regarding eigenvalues of the normalized Laplacian. 

\begin{cor}\label{c:cov}
    Let $G$ and $\tilde{G}$ be weighted graphs with $\tilde{G}$ being a covering of $G$. If $G$ is regular, then the eigenvalues of the normalized Laplacian of $G$ are eigenvalues of the normalized Laplacian of $\tilde{G}$.
\end{cor}

Let $\tilde{B} = (V(\tilde{B}), E(\tilde{B}))$ be a weighted graph with four vertices, $V(\tilde{B}) = \{ v_1, v_2,$ $ v_3, v_4 \}$. The weights of the edges of $\tilde{B}$ are $w(v_1,v_1) = w(v_2,v_2) = n-1$, $w(v_3,v_3)=2$, $w(v_1,v_3)=w(v_2,v_3)=1$, $w(v_3,v_4)=2(n-2)$, $w(v_4,v_4)=2(n-2)(n-1)$, and all other edge weights being zero. This weighted graph is illustrated in Figure~\ref{fig:btilde}. The covering map $p:V(\bpn) \to V(\tilde{B})$ can then be described by the fibers of each the $v_i$, for $i \in \{1,2,3,4\}$, 
    \begin{align*}
        F_1 :=& p^{-1}(v_1) = \{ u \in \bpn :  u(n) = n \},\\
        F_2 :=& p^{-1}(v_2) = \{ u \in \bpn :  u(n) = \overline{n} \},\\
        F_3 :=& p^{-1}(v_3) = \{ u \in \bpn : |u(1)| = n \},\text{ and }\\
        F_4 :=& p^{-1}(v_4) = \{ u \in \bpn : |u(i)| = n, i \not\in \{1,n\} \}.
    \end{align*}

\begin{figure}
    \centering
        \definecolor{rvwvcq}{rgb}{0.08235294117647059,0.396078431372549,0.7529411764705882}
        \begin{tikzpicture}[line cap=round,line join=round,>=triangle 45,scale=0.7]
        
        \draw [line width=1pt] (-5,2)-- (-2,1);
        \draw [line width=1pt] (-5,0)-- (-2,1);
        \draw [line width=1pt] (-2,1)-- (1,1);
        \draw [line width=1pt] (-5.5,2) circle (0.5cm);
        \draw [line width=1pt] (-5.5,0) circle (0.5cm);
        \draw [line width=1pt] (-2,1.5) circle (0.5cm);
        \draw [line width=1pt] (1.5,1) circle (0.5cm);
        \begin{scriptsize}
        \draw [fill=rvwvcq] (-5,2) circle (2.5pt);
        \draw[color=rvwvcq] (-5.35,2) node {$v_1$};
        \draw [fill=rvwvcq] (-5,0) circle (2.5pt);
        \draw[color=rvwvcq] (-5.35,0) node {$v_2$};
        \draw [fill=rvwvcq] (-2,1) circle (2.5pt);
        \draw[color=rvwvcq] (-2,1.35) node {$v_3$};
        \draw [fill=rvwvcq] (1,1) circle (2.5pt);
        \draw[color=rvwvcq] (1.35,1) node {$v_4$};
        \draw[color=black] (-3.53,2.14) node {$1$};
        \draw[color=black] (-3.53,-0.14) node {$1$};
        \draw[color=black] (-0.27,0.44) node {$2(n-2)$};
        \draw[color=black] (-1.25,1.85) node {$2$};
        \draw [fill=rvwvcq] (-5,2) circle (2.5pt);
        \draw[color=black] (-6.65,2) node {$n-1$};
        \draw[color=black] (-6.65,0) node {$n-1$};
        \draw[color=black] (3.65,1) node {$2(n-2)(n-1)$};
        \end{scriptsize}
        \end{tikzpicture}
    \caption{The graph $\tilde{B}$, a projection of $\bpn$}
    \label{fig:btilde}
\end{figure}

We leave it to the reader to verify that the $\tilde{B}$ is in fact a projection, based on Definition~\ref{d:cover}. However, we present the matrices $D$, $A_G$, and $\mathcal{L}$ of $\tilde{B}$.
\begin{equation*}
D = \begin{bmatrix}
    n & 0 & 0 & 0\\
    0 & n & 0 & 0\\
    0 & 0 & 2n & 0\\
    0 & 0 & 0 & 2n(n-2)
\end{bmatrix}, \hfill
A_G = \begin{bmatrix}
    n-1 & 0 & 1 & 0\\
    0 & n-1 & 1 & 0\\
    1 & 1 & 2 & 2(n-2)\\
    0 & 0 & 2(n-2) & 2(n-2)(n-1)
\end{bmatrix}, 
\end{equation*}
\begin{equation*}
\text{ and }\mathcal{L} = \begin{bmatrix}
    \frac{1}{n} & 0 & -\frac{1}{\sqrt{2}n} & 0\\
    0 & \frac{1}{n} & -\frac{1}{\sqrt{2}n} & 0\\
    -\frac{1}{\sqrt{2}n} & -\frac{1}{\sqrt{2}n} & \frac{n-1}{n} & -\frac{\sqrt{n-2}}{n}\\
    0 & 0 & -\frac{\sqrt{n-2}}{n} & \frac{1}{n}
\end{bmatrix}.
\end{equation*}

A direct computation provides the eigenvalues of $\mathcal{L}$, which by Corollary~\ref{c:cov} are eigenvalues of the normalized Laplacian of $\bpn$. Finally, owing to $\bpn$ being regular of degree $n$, we have the following result.
\begin{prop}
    The set of eigenvalues of the normalized Laplacian of $\bpn$ contains $\{0, 1/n, 1/n, 1\}$. Thus, the adjacency spectrum of $\bpn$ contains the values $\{n, n-1, n-1, 0\}$. In particular, the multiplicity of $n-1$ in the spectrum of $\bpn$ is at least 2. 
\end{prop}

Beyond attempting other graph covers/projections we algorithmically generated adjacency matrices of $\bpn$, for small values of $n$, and approximated their eigenvalues. These approximations led to the following conjectures, one of which is a strengthening of our main result. 

\begin{conj}\label{con:one}
    Recall, $\Sp(\bpn)$ is the adjacency spectrum of the graph $\bpn$. Then, $[\ol{(n-1)},n] \subset \Sp(\bpn)$, for all $n\geq3$.
\end{conj}

Our main result establishes part of Conjecture~\ref{con:one}. Moreover, computer evidence suggests the following conjecture. 

\begin{conj}
    Let $\mathop{sg}(\bpn)$ denote the spectral gap of $\bpn$. Then, 
    \begin{enumerate}[(i)]
        \item $\mathop{sg}(\bpn)< 1$, for $n>1$, and
        \item $\mathop{sg}(\bpn)\to 1$ as $n\to\infty$.
    \end{enumerate} 
\end{conj} 

\section*{Acknowledgement} We are indebted to the anonymous referee for their valuable feedback that helped to improve the presentation of this article. 

\bibliographystyle{elsarticle-harv}

\end{document}